\documentclass[final]{siamltex}
\pdfoutput=1
\usepackage{graphicx,bbm}
\usepackage{color}

\def\new#1{#1}

\newtheorem{remark}{Remark}[section]

\def\eps{\varepsilon}

\def\C{\mathbbm{C}}

\def\R{\mathbbm{R}}

\def\Cn{\C^n}

\def\Cnn{\C^{n\times n}}

\def\Re{{\rm Re}}

\def\rank{{\rm rank}}

\def\tr{{\rm tr}}

\def\psa{\sigma_{\kern-.5pt \eps}}
\def\nr{{W}\kern-.75pt}
\def\nabs{\omega}

\def\bmatrix#1{\left[\matrix{#1}\right]}

\mathcode`\@="8000
{\catcode`\@=\active \gdef@{\mkern1mu}} 
\def\BA{{\bf A}}  
\def\Bb{{\bf b}}\def\BB{{\bf B}}

\def\Be{{\bf e}}\def\BE{{\bf E}}  
  
\def\BG{{\bf G}}  
  
\def\BI{{\bf I}}  
  
\def\BK{{\bf K}}

\def\BS{{\bf S}}

\def\BV{{\bf V}}  
  
\def\Bx{{\bf x}}\def\BX{{\bf X}}

\def\BAs{\BA{\kern-1.5pt}}
\def\BLambda{\mbox{\boldmath$\Lambda$}}

\def\mydate{\number\day\ {\ifcase\month \or January\or February\or
              March\or April\or May\or June\or July\or August\or
              September\or October\or November\or December\fi} \number\year}

\title{Fast singular value decay for Lyapunov solutions with nonnormal coefficients}
\author{Jonathan Baker\footnotemark[1]
\and Mark Embree\footnotemark[2]
\and John Sabino\footnotemark[3]}

\begin{document}
\maketitle
\renewcommand{\thefootnote}{\fnsymbol{footnote}}
\footnotetext[1]{Department of Computational and Applied Mathematics, Rice University, 6100 Main Street -- MS 134, Houston, Texas, 77005--1892 ({\tt jpb7@rice.edu}).}
\footnotetext[2]{Department of Mathematics, Virginia Tech, 225 Stanger Street 0123, Blacksburg, Virginia, 24061 ({\tt embree@vt.edu}).}
\footnotetext[3]{The Boeing Company, 9725 E. Marginal Way S, Tukwila, Washington, 98108--4040\\ ({\tt john-paul.n.sabino@boeing.com}).}

\centerline{29 January 2015}\vspace*{5pt}

\begin{abstract}
Lyapunov equations with low-rank right-hand sides often have solutions whose
singular values decay rapidly, enabling 
iterative methods that produce low-rank approximate solutions.
All previously known bounds on this decay involve quantities that depend 
quadratically on the departure of the coefficient matrix from normality:
these bounds suggest that the larger the departure from normality, 
the slower the singular values will decay.
We show this is only true up to a threshold, beyond which
a larger departure from normality can actually correspond
to \emph{faster} decay of singular values:  if the singular values
decay slowly, the numerical range cannot extend far into the
right-half plane.
\end{abstract}

\begin{keywords}
Lyapunov equation, singular values, numerical range, nonnormality
\end{keywords}

\begin{AMS}
15A18, 15A60, 65F30, 93B05, 93B40
\end{AMS}

\pagestyle{myheadings}
\thispagestyle{plain}
\markboth{J. BAKER, M. EMBREE, J. SABINO}
         {FAST SINGULAR VALUE DECAY FOR LYAPUNOV SOLUTIONS}

\section{Introduction} \label{sec:intro}
Lyapunov equations of the form
\begin{equation} \label{eq:lyap}
 \BA\BX + \BX \BAs^* = -\BB\BB^*
\end{equation}
arise from the study of the controllability
and observability of linear time-invariant dynamical systems,
and subsequently in balanced truncation model order reduction~\cite{Ant05b,Zho95}.
In this setting, the right-hand side $-\BB\BB^*$ often has low rank
(equal to the number of inputs or outputs in the system).
If the eigenvalues of $\BA \in \Cnn$ are in the left half of the complex
plane and $(\BA,\BB)$ is controllable,
the solution $\BX\in\Cnn$ is Hermitian positive definite, i.e., $\rank(\BX)=n$~\cite[\S3.8]{Zho95}.  
Even when the coefficient matrix $\BA$ is sparse, $\BX$ is typically dense:
so for large-scale problems one cannot afford to store all $n^2$
entries of the solution.

Penzl observed that, when the right-hand side of~(\ref{eq:lyap}) has low rank,
the singular values $s_1 \ge s_2 \ge \cdots \ge s_n > 0$ of $\BX$ 
often decay exponentially~\cite{Pen00a}, e.g., 
$s_k/s_1 \le C \gamma^k$ for some constants $C>0$ and $\gamma \in (0,1)$.
This fact now enables numerous iterative methods that
seek accurate low-rank approximations to $\BX$; 
\new{see~\cite{BS13,Sim} for recent surveys.}\ \   
Since the singular values of $\BX$ bound the best possible performance
of iterative methods for solving Lyapunov equations, it is important to
understand how they vary with the coefficient matrix $\BA$.\ \ 
Of course, since $\BX$ is Hermitian positive definite, its singular values
equal its eigenvalues; it is common to refer to singular values because
(a) we seek low-rank approximations to $\BX$, and
(b) much of the related analysis generalizes to Sylvester equations, 
where $\BX$ need not even be square.  
We shall thus always speak of the \emph{singular values of $\BX$}, 
$s_1\ge s_2 \ge \cdots \ge s_n > 0$, and the
\emph{eigenvalues of $\BA$}, $\lambda_1,\ldots, \lambda_n\in\sigma(\BA)$.
Let $\C^-$ and $\C^+$ denote the open left and right halves of the complex
plane, and $\|\cdot\|$ denote the vector 2-norm and the matrix norm it 
induces.  
We assume $\BA$ is stable, i.e., $\sigma(\BA)\subset \C^-$.

The singular values of $\BX$ depend on spectral 
properties of $\BA$; grossly speaking, they decay more rapidly 
the farther $\sigma(\BA)$ falls in the left half of the complex plane, 
and more slowly as eigenvalues of $\BA$ grow in imaginary part.  
But eigenvalues alone cannot explain the singular values of $\BX$.
Penzl showed that for any desired singular values of $\BX$, one can
construct a 
corresponding $\BA$ 
with \emph{any spectrum} in the left half-plane
(for some special choice of $\BB$)~\cite{Pen00b}.
Now suppose $\sigma(\BA)$ is fixed.
Recall that $\BA$ is \emph{normal} if it commutes with its adjoint ($\BA\BAs^*=\BAs^*\BA$),
or, equivalently, if eigenvectors give an orthonormal basis for $\C^n$.
We shall use the term \emph{departure from normality} generically; many different
scalar measures of nonnormality have been shown to be essentially equivalent~\cite{EP87}.
All previously known bounds suggest that the singular values of $\BX$ will decay more slowly
as the departure of $\BA$ from normality increases, and it is this 
particular point that concerns us here.  
In Section~\ref{sec:survey} we describe the variety of bounds that 
have been proposed in the literature,
highlighting how they treat the nonnormality of $\BA$.
Section~\ref{sec:example} gives a simple $2\times 2$ example that
clearly illustrates that, in contrast to previously known bounds,
beyond a certain threshold a larger departure from normality can
actually give singular values that \emph{decay more quickly}.
We offer an \new{intuitive} explanation for this behavior in Section~\ref{sec:krylov},
then prove a decay bound that incorporates this effect in Section~\ref{sec:bound}:
the trailing singular values must be small if eigenvalues of the Hermitian 
part of $\BA$ fall far in the right half-plane.

\section{Decay bounds and their inadequacy for nonnormal coefficients} \label{sec:survey}

One approach to proving the decay of the singular values of $\BX$ 
uses the low-rank approximations constructed by the ADI~algorithm; see, e.g.,
\cite{Ant05b,EW86,ES}. 
Suppose $\rank(\BB) = r$.  
The $k$th ADI iteration gives an
approximate solution $\BX_k$ with $\rank(\BX_k)\le k@r$ that satisfies
\[ \BX-\BX_k = \phi_k(\BA)@ \BX@ \phi_k(\BA)^*,\]
where 
\[ \phi_k(z) =  \prod_{j=1}^k {z+\mu_j \over z-\overline{\mu_j}} \]
is a rational function whose parameters, the \emph{shifts} $\{\mu_j\}\subset \C^{+}$, 
are picked from the right half-plane to minimize $\|\BX-\BX_k\|$.
By the optimality of the singular values (the Schmidt--Eckart--Young--Mirsky theorem~\cite[Thm.~3.6]{Ant05b}),
\begin{equation} \label{eq:decay0}
 {s_{kr+1} \over s_1} \le {\| \BX-\BX_k\| \over \|\BX\|} 
                                  \le \|\phi_k(\BA)\|@\|\phi_k(\BAs)^*\|
                                     = \|\phi_k(\BA)\|^2.
\end{equation}
Bounds on the singular values of $\BX$ then follow by approximating 
norms of functions of $\BA$.
Any specific choice of rational function $\phi_k$ gives an 
upper bound, and much theoretical and practical work has addressed
the selection of optimal $\{\mu_j\}$ parameters.  
Since our main point does not depend
on the choice of $\phi_k$, we shall not dwell on that issue here.
Our goal is to illustrate that all known bounds on the singular values 
of $\BX$ fail to capture the diverse behavior possible for nonnormal $\BA$,
so we shall briefly describe the different approaches taken in the literature.
If $\BA$ is normal, then
\begin{equation} \label{eq:fnormal}
 \|\phi_k(\BA)\| = \max_{\lambda\in\sigma(\BA)} |\phi_k(\lambda)|,
\end{equation}
but for nonnormal $\BA$, the left-hand side of~(\ref{eq:fnormal}) 
can be considerably larger than the right-hand side.
There are three common ways to bound $\|\phi_k(\BA)\|$
(cf.~\cite[\S4.11]{Hig08}),
each of which then leads to an upper bound on~(\ref{eq:decay0}). 
\begin{itemize}
\item Eigenvalues:  If $\BA$ is diagonalizable, $\BA=\BV\BLambda\BV^{-1}$, then
\begin{equation} \label{eq:fAeig}
    \|\phi_k(\BA)\| \le \|\BV\|@\|\BV^{-1}\| \max_{\lambda\in\sigma(\BA)} |\phi_k(\lambda)|. 
\end{equation}

Combining~(\ref{eq:fAeig}) with~(\ref{eq:decay0}) gives \\[-9pt]
\begin{equation} \label{eq:szbnd}
 {s_{kr+1}\over s_1} \le \|\BV\|^2\|\BV^{-1}\|^2 
                        \max_{\lambda\in\sigma(\BA)}
                       \prod_{j=1}^k {|\lambda + \mu_j|^2 \over |\lambda - \overline{\mu_j}|^2}.
\end{equation}
This bound was first written down for general diagonalizable $\BA$ by 
Sorensen and Zhou~\cite[Thm.~2.1]{SZ02}, based on earlier work on the
Hermitian case by Penzl~\cite{Pen00b}.  
In that Hermitian case, several concrete bounds have been obtained
by selecting particular real shifts, $\{\mu_j\}$:
using suboptimal shifts, Penzl gave an elegant bound~\cite[Thm.~1]{Pen00b},
which was improved using optimal shifts for a real interval in~\cite[Thm.~2.1.1]{Sab06}.

\ \ \ When $\BA$ is non-Hermitian and the eigenvector matrix is ill-conditioned, 
$\|\BV\| \|\BV^{-1}\|\gg 1$,
one might improve upon~(\ref{eq:szbnd}) 
by posing the maximization problem on larger
subsets of $\C$ that permit constants smaller than $\|\BV\| \|\BV^{-1}\|$.\ \ 
We consider two such methods next.
\item Numerical range: If $\phi_k$ is analytic on the numerical range~\cite[Ch.~1]{HJ91}
\[ \nr(\BA) := \{ \Bx^*\BA\Bx: \Bx\in\Cn, \|\Bx\|=1\},\] 
then 
\begin{equation} \label{eq:fAnr}
  \|\phi_k(\BA)\| \le C \max_{z\in \nr(\BA)} |\phi_k(z)|,
\end{equation}
where \new{$C$ denotes Crouzeix's constant,} $C \in [2,11.08]$~\cite{Cro07}.
Combining this bound with~(\ref{eq:decay0}) gives
\begin{equation} \label{eq:nrbnd}
 {s_{kr+1}\over s_1} \le C^2  
                        \max_{z\in \nr(\BA)}
                       \prod_{j=1}^k {|z+ \mu_j|^2 \over |z - \overline{\mu_j}|^2}.
\end{equation}
This bound only holds when $\phi_k$ is analytic on $\nr(\BA)$, so, in particular,
$\mu_j \not \in \nr(\BA)$.  Since $\mu_j \in \C^+$, a sufficient condition to ensure
analyticity is that $\nr(\BA)\subseteq \C^-$. 
The rightmost extent of $\nr(\BA)$ in the complex plane plays an important
role in analysis of dynamical systems.  
This value is called the \emph{numerical abscissa}
\begin{equation} \label{eq:nabs}
 \nabs(\BA) = \max_{z\in \nr(\BA)} \Re\,z,
\end{equation}
and it equals the rightmost eigenvalue of the Hermitian part of $\BA$:
\[ \nabs(\BA) = \kern-2pt \max_{z\in \nr(\BA)} \kern-2pt {z+\overline{z} \over 2}
              = \max_{\stackrel{\scriptstyle{\Bx\in\Cn}}
                            {\scriptstyle{\|\Bx\|=1}}}
                      \Bx^*\kern-1pt \Big(\kern-.5pt{\BA\kern-.75pt+\kern-.75pt\BAs^*\over 2}\Big) \Bx
              = \max\Big\{ \lambda: \lambda \in \sigma\Big({\BA\kern-.75pt+\kern-.75pt\BAs^*\over 2}\Big)\Big\}.\]
Notice that $\nabs(\BA)$ can be positive even when the spectrum of $\BA$ is in the left half-plane,
and that $|\nabs(\BA)| \le \|\BA\|$.  
The numerical abscissa describes the small $t$ behavior of 
 $\dot\Bx(t) = \BA\Bx(t)$ with $\Bx(0)=\Bx_0$:
\[ \max_{\stackrel{\scriptstyle{\Bx_0\in\Cn}}
                   {\scriptstyle{\|\Bx_0\|=1}}}
        {{\rm d}\over {\rm d} t} \|\Bx(t)\|\Big|_{t=0} = \nabs(\BA);\]
see, e.g., \cite[Thm.~17.4]{TE05}.  Thus $\nabs(\BA)>0$ is a necessary condition
for solutions of $\dot{\Bx}(t) = \BA\Bx(t)$ to exhibit transient growth.

\item Pseudospectra:
\new{The requirement in~(\ref{eq:nrbnd}) that $\phi_k$ be 
analytic throughout $W(\BA)$, which is typically reduced to 
the condition $\nr(\BA)\subset \C^-$, excludes many stable $\BA$.  
Thus we consider a more flexible alternative.}
Given $\eps>0$, if $\phi_k$ is analytic on the $\eps$-pseudsopectrum~\cite{TE05}
\begin{eqnarray*}
 \psa(\BA) &=& \{ z\in\C: \mbox{$z\in \sigma(\BA)$ or $\|(z-\BA)^{-1}\|>1/\eps$}\}\\
&=& \{ z\in\C: \mbox{$z\in\sigma(\BA+\BE)$ for some $\BE\in\Cnn$ with $\|\BE\|<\eps$}\},
\end{eqnarray*}
then 
\begin{equation} \label{eq:fApsa}
    \|\phi_k(\BA)\| \le {L_\eps \over 2 \pi \eps} \sup_{z\in\psa(\BA)} |\phi_k(z)|,
\end{equation}
where $L_\eps$ denotes the contour length of the boundary of $\psa(\BA)$; see, 
e.g., \cite[p.~139]{TE05}.
Substituting~(\ref{eq:fApsa}) into~(\ref{eq:decay0}) yields~\cite[(3.4)]{Sab06}
\begin{equation} \label{eq:penzl}
 {s_{kr+1}\over s_1} \le {L_\eps^2 \over 4 \pi^2 \eps^2}
                        \max_{z\in \psa(\BA)}
                       \prod_{j=1}^k {|z+ \mu_j|^2 \over |z - \overline{\mu_j}|^2}.
\end{equation}
The choice of $\eps>0$ balances the leading constant
against the set over which the maximization occurs:  increasing $\eps$ typically
decreases $L_\eps^2/(4\pi^2\eps^2)$ but enlarges $\psa(\BA)$.
For any $\{\mu_j\}\subset \C^+$ there exists $\eps>0$ sufficiently
small that $\phi_k$ is analytic on $\psa(\BA)$, since $\sigma(\BA)\subset \C^-$
and $\psa(\BA)$ converges to $\sigma(\BA)$ in the Hausdorff metric as $\eps\to 0$;
\new{see, e.g., \cite[Ch.~4]{TE05} for related details.}
\end{itemize}

All these bounds derived from~(\ref{eq:decay0}) predict the decay of singular
values will slow as the departure of $\BA$ from normality increases,
as reflected in increased ill-conditioning of the eigenvector matrix
(i.e., the eigenvectors associated with distinct eigenvalues become increasingly
aligned), enlargement of the numerical range, or an increase in the sensitivity
of the eigenvalues to perturbations.

Several alternative bounds on the singular values of $\BX$ 
have been derived using entirely different approaches,
but they share this same property.
\new{
Convergence theorems for the low-rank approximate solutions to
the Lyapunov equation
constructed by projection methods also provide upper bounds on
the decay of the singular values of $\BX$.\ \  In this literature,
results based on the numerical range have proved most popular.
Like~(\ref{eq:nrbnd}), these bounds predict slower singular value
decay as the distance of the numerical range from the origin decreases, 
and they fail to hold when $0\in\nr(\BA)$ (see, e.g., 
Theorem~4.2 of~\cite{DKS11}, which resembles~(\ref{eq:nrbnd}),
and Corollary~2.5 of~\cite{Bec11}). Thus they do not apply
to the highly nonnormal examples that interest us here.}

\new{Antoulas, Sorensen, and Zhou~\cite[Thm.~3.1]{ASZ02} propose a different
strategy.  For diagonalizable $\BA$ 
they write the solution $\BX$ as a finite series to show, for $r=1$,}
\begin{equation} \label{eq:asz}
  s_{k+1} \le (n-k)^2 \|\BV\|^2\|\BV^{-1}\|^2 \|\BB\|^2\, \delta_{k+1}, 
\end{equation}
where 
\[ \delta_{k} = -{1\over 2\,\Re\,\lambda_{k}} \prod_{j=1}^{k-1} 
                   {|\lambda_k-\lambda_j|^2\over |\lambda_k + \overline{\lambda_j}|^2},\]
with the eigenvalues $\lambda_1, \ldots, \lambda_n$ of $\BA$ ordered to make 
$\delta_1 \ge \delta_2 \ge \cdots \ge \delta_n$.
By~(\ref{eq:lyap}), we have
\begin{equation} \label{eq:s1}
 \|\BB\|^2 = \|\BB\BB^*\| = \|\BA\BX+\BX\BA^*\| \le 2\|\BA\| \|\BX\| = 2@ \|\BA\|@ s_1,
\end{equation}
so~(\ref{eq:asz}) 
implies the relative bound
\[ 
  {s_{k+1}\over s_1}  \le 2 (n-k)^2 \|\BA\|@ \|\BV\|^2\|\BV^{-1}\|^2\, \delta_{k+1}.\]
(The $r>1$ case is slightly more complicated~\cite[Thm.~3.2]{ASZ02}.)
By analyzing the trace of $\BX$, Truhar and Veseli\'{c}~\cite{TV07} 
derive an alternative to~(\ref{eq:asz})
that characterizes the departure of $\BA$ from normality by terms like $\|\BV\|^2 \|\widehat{\Bb}_j\|^2$, 
where $\widehat{\Bb}_j^*$ denotes the $j$th row of $\BV^{-1}\BB$ 
(and thus depends on the conditioning of the eigenvectors of $\BA$).  
This bound can be 
generalized to nondiagonalizable $\BA$, with an explicit formula given for
$2\times 2$ Jordan blocks~\cite[Thm.~2.2]{TV07}, and can be further
generalized to Sylvester equations~\cite{TTL10}.  
Bounds for coefficients $\BA$ that are non-self-adjoint operators on Hilbert space 
exhibit similar dependence on the square of the condition number of the transformation 
that orthogonalizes a Riesz basis of eigenvectors~\cite[Thm.~4.1]{GK14}.

\begin{figure}[t!]
\begin{center}
\includegraphics[scale=.37]{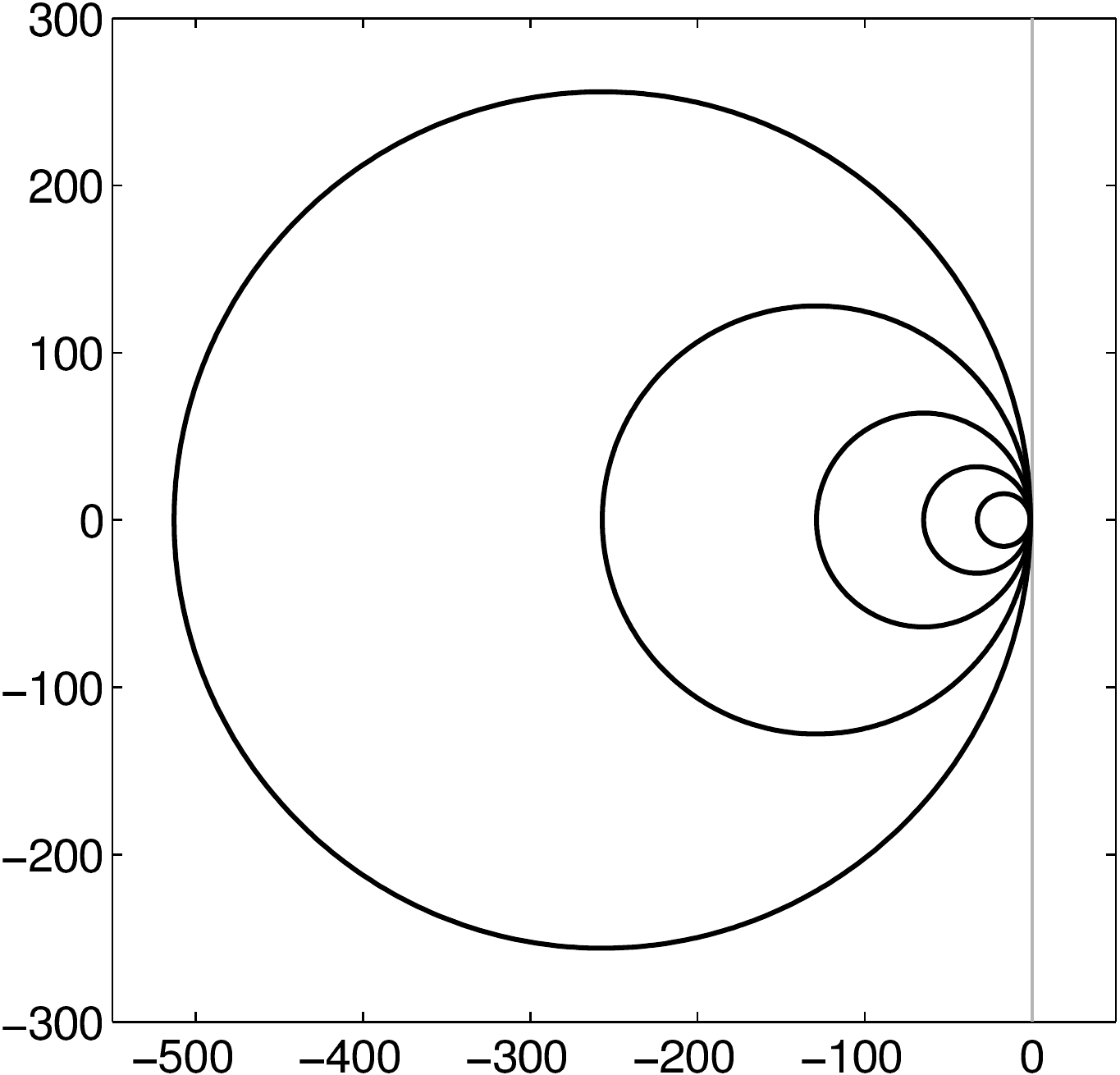}
\qquad\quad
\includegraphics[scale=.37]{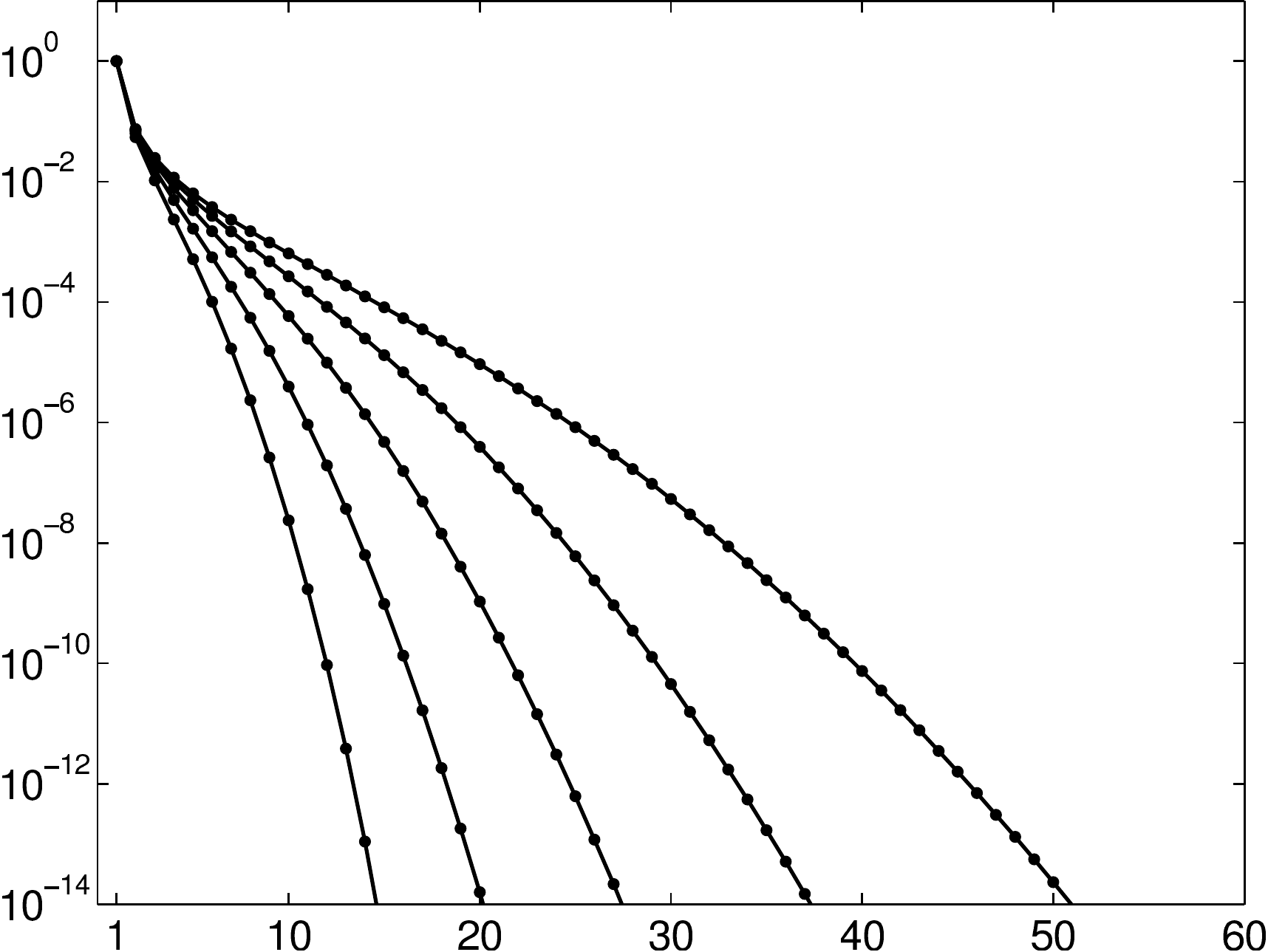}
\begin{picture}(0,0)
\put(-92,-10){\small $k$}
\put(-203,82){\small $\displaystyle{s_k\over s_1}$}
\put(-148,42){\rotatebox{-80}{\small $n=2^4$}}
\put(-57,42){\rotatebox{-49}{\small $n=2^8$}}
\put(-266,15){\rotatebox{0}{\small $n=2^8$}}
\put(-273,36){\rotatebox{0}{\small $n=2^7$}}
\end{picture}
\end{center}
\vspace*{4pt}
\caption{\label{fig:dtest}
Boundaries of $\nr(\BA)$ (left) and decay of singular values of $\BX$ 
(right) for discretizations of ${\rm d}/{\rm d} x - 1$ of dimension 
$n = 2^4, 2^5, \ldots, 2^8$.  As $n$ increases, the numerical range 
enlarges, while the singular values decay more slowly.
This correlation is consistent with previously known bounds.
}
\end{figure}
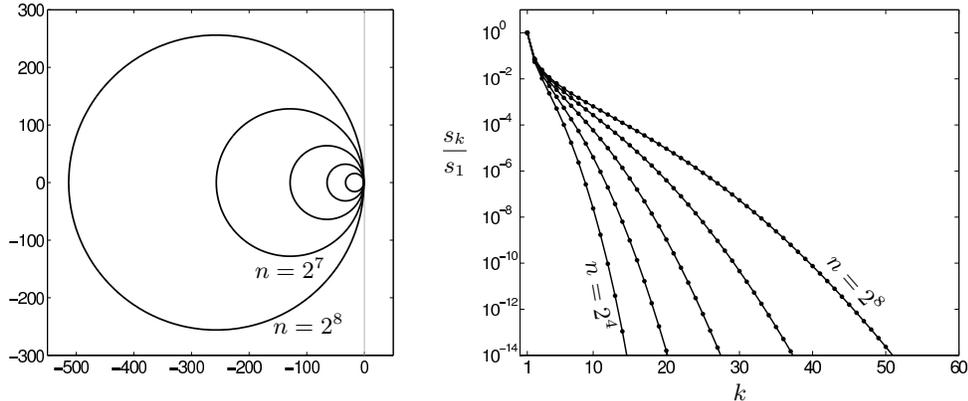

\medskip
When $\nr(\BA)\subset \C^-$, these bounds can be qualitatively descriptive,
even when $\BA$ departs significantly from normality.  
For a simple example, suppose $\BA$ is a discretization of the differential operator
\[ {{\rm d} \over {\rm d} x} - 1\]
defined on absolutely continuous functions in $L^2(0,1)$ satisfying $u(1)=0$.
Approximating the operator with forward finite differences on the uniform grid 
with spacing $1/n$ gives
\[ \BA = \bmatrix{-1-n & n \cr & -1-n & \ddots \cr & & \ddots & n \cr & & & -1-n}
         \in \Cnn\]
with spectrum $\sigma(\BA) = \{-1-n\}$ in the left half-plane.  Since $\BA$ 
is a Jordan block, its numerical range is known in closed form~\cite{MS79}:
\[ \nr(\BA) = \Big\{ z\in\C: |z+1+n| \le n \cos\Big({\pi\over n+1}\Big) \Big\},\]
a disk centered at $-1-n$ of radius \new{$n\cos(\pi/(n+1))$}.
Notice that as $n$ increases $\nr(\BA)$ enlarges monotonically:
the numerical range includes larger portions of the 
half-plane $\{ z\in\C: \Re\,z<-1\}$, reflecting the resolvent behavior
of the underlying differential operator~\cite[\S5]{TE05}.
As $n$ increases, the singular value decay slows.
This behavior is shown in Figure~\ref{fig:dtest}, where
$\BB$ is a constant vector.  
In this case, as predicted by the bounds we have surveyed,  
an increasing departure from normality slows convergence.
We shall see that 
$\nabs(\BA) = -1-n(1-\cos(\pi/(n+1)))<0$ is a crucial property.

\medskip
Not all nonnormal coefficients give this same behavior.  
To see how the known bounds fail to capture the rich behavior exhibited by the 
singular values of Lyapunov solutions with highly nonnormal coefficients,
consider those $\BX$ that exhibit no decay at all, i.e.,
$\BX = \xi@\BI$ for some $\xi>0$,
for ${\rm rank}(\BB)<n$.
Then~(\ref{eq:lyap}) reduces to
\[ \BA + \BAs^* = {- 1\over \xi} \BB\BB^*,\]
which implies that the Hermitian part of $\BA$ is a negative semidefinite
matrix with rightmost eigenvalue (hence numerical abscissa, $\nabs(\BA)$), 
equal to zero.  
If the numerical abscissa is positive, reflecting a larger departure from
normality, the singular values \emph{must} decay faster.
Important applications give rise to matrices with $\nabs(\BA)>0$;
for example, positive $\nabs(\BA)$ can grow with Reynolds number in fluid flows,
a fact that complicates studies of transition to turbulence~\cite{TTRD93}.
Lyapunov equations with low-rank right-hand sides have recently been 
applied to study this problem~\cite{EMSW12}.
To cleanly illustrate the inadequacy of existing bounds, 
we next study a family of $2\times 2$ matrices.
\section{A completely solvable example} \label{sec:example}
Consider the following $2\times 2$ example from~\cite{Sab06}, 
where we interpret ``singular value decay'' to mean the ratio of
the first two singular values, $s_2/s_1$.
Consider the coefficient and right-hand side%
\footnote{Note the normalization of $\BB$; if the second component of $\BB$ is zero,
then $\BB$ is an eigenvector of $\BA$, and the corresponding linear system is not
controllable~\cite{Ant05b}.}
\[ \BA(\alpha) = \bmatrix{-1 & \alpha \cr 0 & -1}, \qquad \BB = \bmatrix{t \cr 1}.\]
Note that $W(\BA(\alpha))$ is the disk in $\C$ centered at $\lambda=-1$ with radius $|\alpha|/2$. 
The solution to the Lyapunov equation can be written out explicitly:
\[ \BX = {1\over 4} \bmatrix{2t^2 + 2\alpha t + \alpha^2 & \alpha + 2 t \cr
                     \alpha+2t & 2}.\]
We seek the the right-hand side $\BB$ that gives the \emph{slowest} decay, i.e.,
that maximizes the ratio
\[ {s_2\over s_1}
    = {\tr(\BX) - \sqrt{\tr(\BX)^2-4\,\det(\BX)}
       \over  \tr(\BX) +\sqrt{\tr(\BX)^2-4\,\det(\BX)}} \le 1\]
over all controllable $\BB\in\R^2$, i.e., over all $t\in\R$.
This worst case decay is attained when $t = -\alpha/2$, giving
\[ {s_2 \over s_1}
         = \left\{ \begin{array}{ll} \alpha^2/4, & 0<\alpha \le 2;\\
                                     4/\alpha^2, & 2\le \alpha.
                   \end{array}\right.\]
As $\alpha$ increases from zero, so too does the departure of $\BA(\alpha)$ from normality.
The ratio $s_2/s_1$ also increases, but only up to $\alpha=2$ (when $\nabs(\BA(\alpha))=0$).
As $\alpha$ increases beyond $\alpha=2$, the ratio of singular values decreases significantly:
contrary to our expectation from bounds described in Section~\ref{sec:survey}, 
\emph{the decay actually improves}.

\section{Krylov conditioning and decay} \label{sec:krylov}

We can gain some general insight into this decay behavior by writing
a Lyapunov solution $\BX$ in terms of the solution of a related
canonical Lyapunov equation that only depends on the spectrum of $\BA$.
\new{This formulation is not intended for practical calculations, but it provides
some intuition for the results that follow in Section~\ref{sec:bound}.}

Let $\BB\in\C^{n\times 1}$ and suppose $(\BA,\BB)$ is controllable.  
Thus $\BA$ is nonderogatory, so its minimum polynomial equals its characteristic polynomial,
\[ \chi(z) = (z-\lambda_1) \cdots (z-\lambda_n) 
           = c_0 + c_1 z + \cdots + c_{n-1}z^{n-1} +z^n.\]
Let $\BA_c$ be the associated companion matrix,
\[ \BA_c = \bmatrix{& & & -c_0 \cr 1 & & & -c_1 \cr & \ddots & & \vdots \cr & & 1 & -c_{n-1}},  \]
whose eigenvalues are the same as those of $\BA$.\ \ Antoulas, 
Sorensen, and Zhou~\cite[Lem.~3.1]{ASZ02} describe the following
method for constructing the solution $\BX$ to $\BA\BX+\BX\BAs^*=-\BB\BB^*$.
Let $\BK$ denote the Krylov matrix
\[ \BK = [\BB\ \BA\BB\ \cdots\ \BAs^{n-1}\BB] \in \C^{n\times n},\]
and $\Be_1$ be the first column of the $n\times n$ identity matrix.
Then $\BA\BX+\BX\BAs^* = -\BB\BB^*$ if and only if $\BX = \BK\BG\BK^*$, where
$\BG$ solves the companion Lyapunov equation
\[ \BA_c^{} \BG + \BG\BA_c^{\kern-1.5pt *} = -\Be_1^{}\Be_1^*.\]
Notice that $\BG$ depends only on $\BA_c$, and hence only on the spectrum of $\BA$,
not the departure of $\BA$ from normality or the right-hand side $\BB$:
the influence of these latter factors on $\BX$ occurs only through the matrix $\BK$.

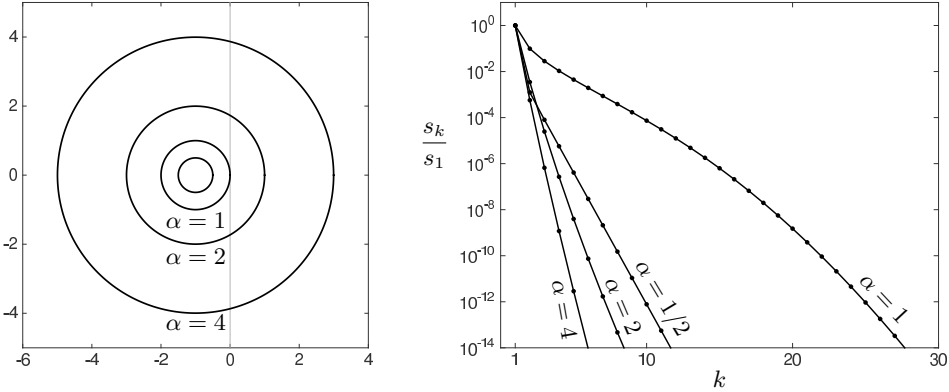
\begin{figure}[t!]
\begin{center}
\includegraphics[scale=.37]{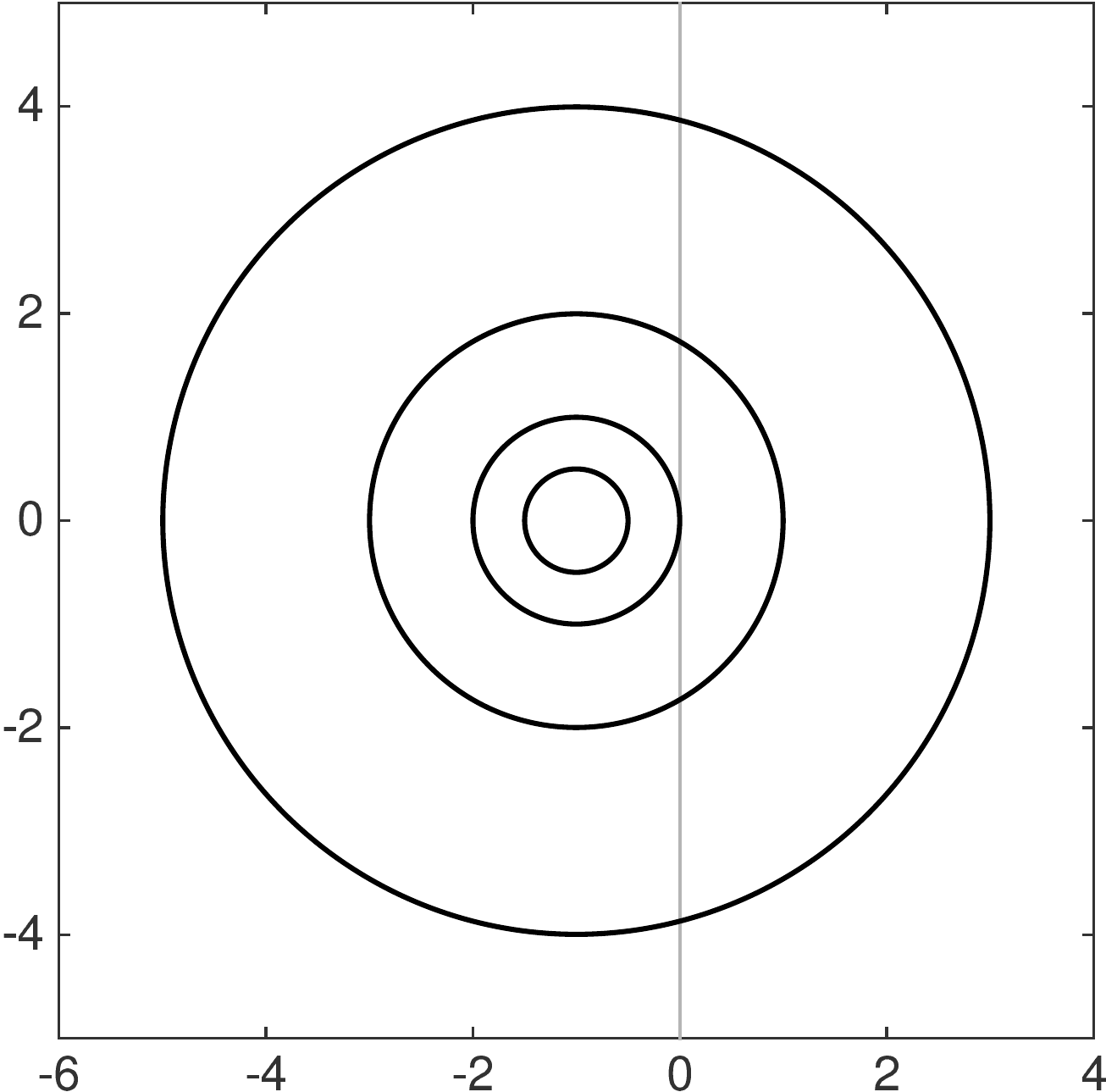}
\qquad\quad
\includegraphics[scale=.37]{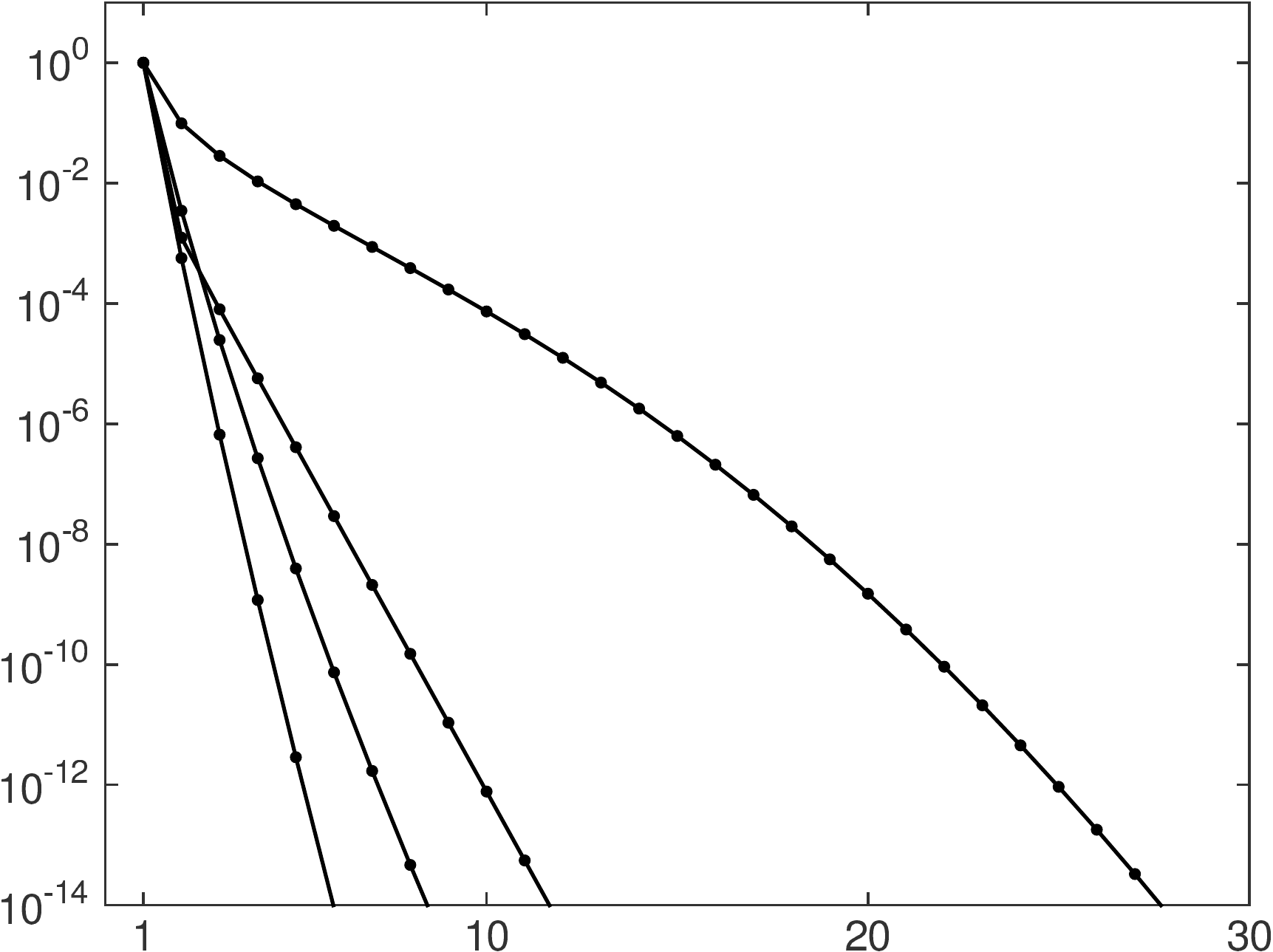}
\begin{picture}(0,0)
\put(-92,-8){\small $k$}
\put(-203,82){\small $\displaystyle{s_k\over s_1}$}
\put(-153,32){\rotatebox{-80}{\small $\alpha=4$}}
\put(-133,32){\rotatebox{-67}{\small $\alpha=2$}}
\put(-123.5,38){\rotatebox{-59}{\small $\alpha=1/2$}}
\put(-37,32){\rotatebox{-47}{\small $\alpha=1$}}
\put(-299,52){\rotatebox{0}{\small $\alpha=1$}}
\put(-299,38.5){\rotatebox{0}{\small $\alpha=2$}}
\put(-299,13.5){\rotatebox{0}{\small $\alpha=4$}}
\end{picture}
\end{center}
\vspace*{3pt}
\caption{\label{fig:atest}
Boundaries of $\nr(\BA)$ (left) and decay of singular values of $\BX$ (right)
for Jordan blocks~(\ref{eq:jblock}) of dimension $n=64$ with off-diagonal
$\alpha = 1/2, 1, 2, 4$.  As $\alpha$ increases, the numerical range enlarges.  
For small and large $\alpha$, the singular values of $\BX$ decay quickly;
the $\alpha=1$ case, having an intermediate departure of $\BA$ from normality,
gives singular values of $\BX$ that decay much slower. (Here $\BB = [1, \ldots, 1]^*$.)}
\vspace*{-10pt}
\end{figure}

Let $\varsigma_k(\cdot)$ denote the $k$th singular value of a matrix.
Since $\BG$ is positive definite, it has a square root, and so
\[ s_k := \varsigma_k(\BX) = \varsigma_k(\BK\BG\BK^*) = \varsigma_k(\BK\BG^{1/2})^2 
                    \le \varsigma_k(\BK)^2 \varsigma_1(\BG^{1/2})^2 
                     = \varsigma_k(\BK)^2 \|\BG\|,\]
using the singular value inequality~\cite[Thm.~3.3.16(d)]{HJ91}.
Use~(\ref{eq:s1}), $\|\BB\BB^*\| \le 2 \|\BA\|@s_1$, 
to obtain the bound
\begin{equation} \label{eq:krybnd}
   {s_k\over s_1} \le  \varsigma_k(\BK)^2 \|\BA\| \bigg({2@\|\BG\|\over \|\BB\BB^*\|}\bigg).
\end{equation}
The singular values of $\BX$ will thus decay (at least) at a rate controlled by 
the singular values of the Krylov matrix $\BK$;
note that the term in parentheses in~(\ref{eq:krybnd}) is independent of the departure
of $\BA$ from normality.
The columns of $\BK$ are iterates of the power method, hence one can gain insight into the
decay of singular values of $\BX$ by studying the convergence of the power method 
for nonnormal~$\BA$.
(See~\cite[\S28]{TE05}, especially the illustration in Fig.~28.1 showing 
how nonnormality can \emph{accelerate} the convergence of the power method.)
We shall not pursue this direction here, but instead imagine fixing $\BB$ and the 
spectrum of $\BA$, then varying the departure of $\BA$ from normality, 
e.g., 
\[ \BA = \BLambda + \alpha \BS,\]
where $\BLambda$ is diagonal, $\BS$ is strictly upper triangular, and $\alpha$ controls
the departure of $\BA$ from normality.
For a concrete example, take $\BLambda=-\BI$ and $\BS$ to be the shift matrix,
yielding a Jordan block that generalizes the example in Section~\ref{sec:example}:
\begin{equation} \label{eq:jblock}
 \BA = \bmatrix{-1 & \alpha \cr  & -1 & \ddots \cr & & \ddots & \alpha \cr & & & -1}.
\end{equation}
The departure of $\BA$ from normality is small when $\alpha$ is small.
In this case $\BA \approx -\BI$, so all columns of $\BK$ will be 
nearly the same: $\varsigma_k(\BK)$ will be small for all $k\ge 2$,
and~(\ref{eq:krybnd}) captures the fast decay of the singular values of $\BX$.
For large $\alpha$, the matrix $\BK$ will be severely graded\new{; specifically,} the
norm of each column of $\BK$ will be on the order of $\alpha^{k-1}$.
Thus for large $\alpha$, the singular values of $\BK$ must also decay
rapidly.%
\footnote{One could apply results on the singular values of graded matrices, e.g., \cite{Ste01},
to obtain quantitative estimates.}
\new{Since $\|\BA\|$ only grows linearly with $\alpha$,%
\footnote{\new{Gerschgorin's theorem applied to $\BAs^*\BA$ gives
$\alpha-1\le\|\BA\|\le \alpha+1$ for $\alpha> 3$.}}
by~(\ref{eq:krybnd}), the singular values
of $\BX$ must decay quickly as well.}
The slowest decay should thus occur for values of $\alpha$ that are neither 
too small nor too large, as suggested by the two dimensional case.
Indeed, this intuition is confirmed in Figure~\ref{fig:atest}, which
shows an example with $n=64$ and $\alpha=1/2, 1, 2, 4$.  Of the 
cases shown, the singular values decay most slowly for $\alpha=1$, 
when the rightmost extent of $W(\BA)$ comes closest to the
imaginary axis.
We next describe rigorous bounds that connect properties of $W(\BA)$
to the decay of the singular values of $\BX$.

\section{Large numerical abscissa implies fast decay} \label{sec:bound}

In~(\ref{eq:nabs}) we defined the numerical abscissa, $\omega(\BA)$, 
which is both the rightmost extent of the numerical range and the rightmost
eigenvalue of the Hermitian part $(\BA+\BAs^*)/2$ of $\BA$.  
The subordinate eigenvalues of the Hermitian part further
inform our understanding of the departure of $\BA$ from normality.  
For example, these eigenvalues have recently been used to
bound the number of Ritz values of $\BA$ that can fall in subregions
of $\nr(\BA)$~\cite[Thm.~1.2]{CE12}.%
\footnote{In the context of moment-matching model reduction algorithms~\cite[Ch.~11]{Ant05b},
these results relating Ritz values to the eigenvalues of $(\BA+\BAs^*)/2$ restrict the number of 
poles of a reduced-order model that can fall in the right half-plane.}
Like $\omega(\BA)$, interior eigenvalues of $(\BA+\BAs^*)/2$ 
can be positive even when $\BA$ is stable.  
The following theorem bounds these eigenvalues in terms of the 
singular values of $\BX$.\ \ 
This result can be read from two different perspectives:
given the singular values of $\BX$, the bound reveals 
something about those $\BA$ that can support such solutions
(Theorem~\ref{thm:genbnd} and Corollary~\ref{cora:s1n});
given $\BA$, one obtains an upper bound on the decay of singular
values of $\BX$ that requires, in a specific context, 
\emph{faster} decay as the departure \new{of $\BA$ from normality increases}
(Corollary~\ref{cora:genbnd}).

\medskip
\begin{theorem} \label{thm:genbnd}
Let $\BX\in\Cnn$ solve the Lyapunov equation~$(\ref{eq:lyap})$ with $(\BA,\BB)$ 
controllable.  Then for all $k=1,\ldots, n$,
\begin{equation} \label{eq:genbndom}
     {s_k\over s_1} - 1 -{\|\BB\|^2 \over 2 s_1 \|\BA\|} \le {\omega_k \over \|\BA\|} \le 1 - {s_{n-k+1} \over s_1},
\end{equation}
where $\omega_k$ denotes the $k$th rightmost eigenvalue of ${1\over 2}(\BA+\BAs^*)$
and $s_k$ denotes the $k$th singular value of $\BX$.
\end{theorem}

{\em Proof}.
Write the solution $\BX = \xi(\BI-\BE)$ for $\xi> 0$ and $\BE$ Hermitian.
Then since $\BX$ solves the Lyapunov equation~(\ref{eq:lyap}),  
\begin{equation} \label{eq:HAeig}
 {\BA+\BAs^*\over 2} = -{1\over 2\xi} \BB\BB^*+{\BA\BE+\BE\BAs^* \over 2}.
\end{equation}
Let $\lambda_k(\cdot)$ denote the $k$th eigenvalue of a Hermitian matrix,
labeled from right to left, and let $\varsigma_k(\cdot)$ 
the $kth$ singular value of a matrix, again labeled from largest to smallest.
Weyl's inequalities for the eigenvalues of sums of Hermitian matrices
(see, e.g., \cite[Thm.~4.3.1]{HJ13}) imply
\[ \lambda_n\Big(-{1\over 2\xi} \BB\BB^*\Big) 
   + \lambda_k\Big({\BA\BE+\BE\BAs^* \over 2}\Big)
   \le 
   \lambda_k\Big(-{1\over 2\xi} \BB\BB^*+{\BA\BE+\BE\BAs^* \over 2}\Big) \]
and 
\[ \lambda_k\Big(-{1\over 2\xi} \BB\BB^*+{\BA\BE+\BE\BAs^* \over 2}\Big) 
    \le
     \lambda_1\Big(-{1\over 2\xi} \BB\BB^*\Big) 
   + \lambda_k\Big({\BA\BE+\BE\BAs^* \over 2}\Big).
\]
Since $-\BB\BB^*/2\xi$ is Hermitian negative semidefinite, 
\[ 
\lambda_n\Big(-{1\over 2\xi} \BB\BB^*\Big) = - {\|\BB\|^2 \over 2\xi}, \qquad
\lambda_1\Big(-{1\over 2\xi} \BB\BB^*\Big) \le 0. 
\]
Now by equation~(\ref{eq:HAeig}),
\[ \lambda_k\Big(-{1\over 2\xi} \BB\BB^*+{\BA\BE+\BE\BAs^* \over 2}\Big)
     = 
 \lambda_k\Big({\BA+\BAs^*\over 2}\Big) =: \omega_k.
\] 
Together, these pieces imply
\begin{equation} \label{eq:omineq}
-{\|\BB\|^2 \over 2 \xi} + \lambda_k\Big({\BA\BE+\BE\BAs^* \over 2}\Big)
 \le \omega_k 
    \le \lambda_k\Big({\BA\BE+\BE\BAs^* \over 2}\Big).
\end{equation}
Note that $(\BA\BE+\BE\BAs^*)/2$ is the Hermitian part of $\BA\BE$.
The $k$th singular value of a matrix gives an upper bound on the
$k$th rightmost eigenvalue of its Hermitian part~\cite[Cor.~3.1.5]{HJ91}.
Applying this bound to both $\BA\BE$ and $-\BA\BE$ gives
\[ -\varsigma_{n-k+1}(\BA\BE) \le  \lambda_k\Big({\BA\BE+\BE\BAs^* \over 2}\Big) \le \varsigma_k(\BA\BE).\]
Using the singular value inequality~\cite[Thm.~3.3.16(d)]{HJ91},
\[ \varsigma_k(\BA\BE)\le \varsigma_1(\BA)\,\varsigma_k(\BE) = \|\BA\|\, \varsigma_k(\BE),\]
obtain from~(\ref{eq:omineq}) that
\begin{equation} \label{eq:omkbnd}
-{\|\BB\|^2 \over 2 \xi \|\BA\|} - \varsigma_{n-k+1}(\BE) 
\le {\omega_k \over \|\BA\|} \le \varsigma_k(\BE).
\end{equation}
Since $\BE = \BI - \BX/\xi$, the eigenvalues of $\BE$, labeled from right to left, are
\[ \lambda_k(\BE) = 1- s_{n-k+1}/\xi, \qquad \rlap{$k=1,\ldots, n$.}\]
The form $\BX=\xi(\BI-\BE)$ allows for various choices of $\xi$ and $\BE$.  
Taking $\xi=s_1$ gives $\BE = \BI-\BX/s_1$, hence
$0 = \lambda_n(\BE)\le \cdots \le \lambda_1(\BE)$ and
\[ \varsigma_k(\BE) = 1-s_{n-k+1}/s_1.\]
Thus~(\ref{eq:omkbnd}) implies
\[ 
     {s_k\over s_1} - 1 -{\|\BB\|^2 \over 2 s_1 \|\BA\|} \le {\omega_k \over \|\BA\|} \le 1 - {s_{n-k+1} \over s_1}. \rlap{\qquad \endproof}
\]

\medskip
\begin{remark} \rm
In the proof of Theorem~\ref{thm:genbnd}, the choice $\xi=s_1$ for the scaling factor
$\xi$ is usually suboptimal.  Smaller values of $\xi>0$ can give tighter bounds but
usually at the expense of more intricate formulas (since then the eigenvalues 
of $\BE$ can be positive and negative).  
As a special case, we can take $\xi = (s_1+s_n)/2$ to 
optimize~(\ref{eq:omkbnd}) for $k=1$,
giving $\lambda_1(\BE) = -\lambda_n(\BE) = (s_1-s_n)/(s_1+s_n)$ 
and
\[ \nabs(\BA) \le {s_1-s_n \over s_1+s_n} \|\BA\|.\]
This expression has a nice interpretation: if the smallest singular value $s_n$ of $\BX$
is on the same order as $s_1$, then $\nabs(\BA)$ must be quite a bit smaller than $\|\BA\|$.
When combined with the $k=1$ lower bound from Theorem~\ref{thm:genbnd} 
(with $\xi=s_1$), we obtain bounds on the rightmost extent of any
numerical range that can support a solution $\BX$ with extreme
singular values $s_1$ and $s_n$.
\end{remark}

\medskip
\begin{corollary} \label{cora:s1n}
For controllable $(\BA,\BB)$, the numerical abscissa $\nabs(\BA)$ is bounded 
by the extreme singular values of the solution $\BX\in\Cnn$ to the Lyapunov 
equation~$(\ref{eq:lyap})$: 
\begin{equation} \label{eq:s1n}
  -{\|\BB\|^2 \over 2s_1} \le \nabs(\BA) \le  {s_1-s_n \over s_1+s_n} \|\BA\|.
\end{equation}
\end{corollary}

\begin{figure}[t!]
\begin{center}
\includegraphics[scale=0.45]{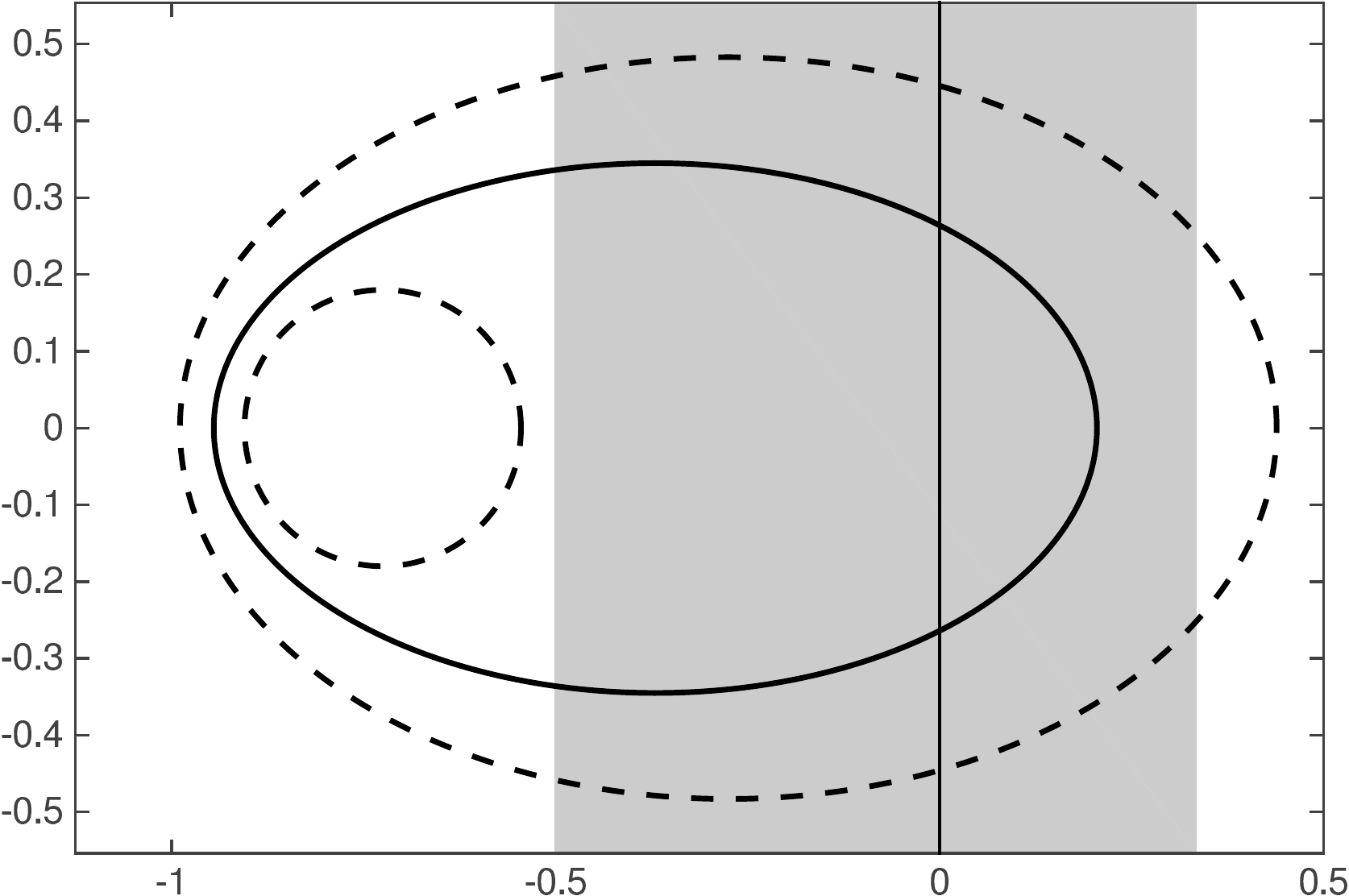}
\end{center}
\vspace*{-.25em}
\caption{\label{fig:strip}
Visualization of Corollary~\ref{cora:s1n} in the complex plane, 
with $\|\BA\|=\|\BB\|=s_1=1$ and $s_n = 1/2$.  
To allow such slow decay of the singular values of $\BX$, by~(\ref{eq:s1n}) 
the rightmost extent of the numerical range, $\nabs(\BA)$,  must fall within 
the gray strip.
\new{The solid curve shows the boundary of one such numerical range.}
The dashed curves show the boundaries
of two different numerical ranges for which the singular values must decay more rapidly,
since in each case the numerical abscissa violates~(\ref{eq:s1n}).
}
\vspace*{-.25em}
\end{figure}

Figure~\ref{fig:strip} provides a schematic illustration of this Corollary.
\new{When the singular values decay slowly (as described in the caption), the
rightmost extent of the numerical range must fall within the gray strip.
Note that the converse need not hold:  the singular values can decay quickly
regardless of $\nr(\BA)$, depending on $\BB$ and finer spectral properties
of $\BA$.}

\medskip
Rearranging the upper bound in Theorem~\ref{thm:genbnd} gives an upper bound 
on the decay of the trailing singular values of $\BX$.

\medskip
\begin{corollary} \label{cora:genbnd}
For controllable $(\BA,\BB)$, 
the singular values of the solution $\BX\in\Cnn$ to the Lyapunov equation~$(\ref{eq:lyap})$ satisfy
\begin{equation} \label{eq:genbnd}
    {s_{n-k+1} \over s_1} \le 1 - {\omega_k \over \|\BA\|},
     \rlap{\qquad $k=1,\ldots, n$.}
\end{equation}
\end{corollary}

\begin{remark}\rm
As observed in Section~\ref{sec:survey}, the case of no decay ($s_1=s_n$)
implies that $\omega_1 \equiv \omega(\BA) = 0$, in which case 
Corollary~\ref{cora:genbnd} with $k=1$ is sharp.
On the other hand, in the highly nonnormal case where $0 < \omega_k \approx \|\BA\|$, 
Corollary~\ref{cora:genbnd} requires that the $k$th lowest singular value be
small, \emph{regardless of $\BB$}.  This stands in contrast to the traditional
bounds surveyed in Section~\ref{sec:survey} for two reasons:
higher nonnormality implies faster decay, rather than slower decay;
the rank of $\BB$ does not feature in the bound on $s_{n-k+1}/s_1$, whereas 
the other bounds predict slower decay as the rank of $\BB$ increases.

Corollary~\ref{cora:genbnd} is designed to show that decay must occur in this 
specific highly nonnormal scenario.
The result is not useful when $\|\BA\|$ is controlled by eigenvalues far in the 
left half-plane, rather than being dominated by the departure of $\BA$ from normality.
In this case $s_{n-k+1}/s_1$ can be quite small while the right-hand side of~(\ref{eq:genbnd})
is not.  In particular, when $\omega_k<0$ (as must occur for all $k$ 
when $\BA$ is stable and normal), the bound in~(\ref{eq:genbnd}) is vacuous.

\end{remark}

\medskip
\begin{remark}\rm
Note that the rate of decay could be even stronger than indicated by
Corollary~\ref{cora:genbnd}.  For the $2\times 2$ Jordan block
considered in Section~\ref{sec:example},
\[ \omega_1 = \alpha/2-1, \qquad 
   \|\BA\| = \sqrt{1 + \alpha^2/2 + \alpha\sqrt{\alpha^2/4 + 1}},\]
so Corollary~\ref{cora:genbnd} gives the bound
\[ {s_2\over s_1} \le 1 - {\omega_1 \over \|\BA\|} \to 
       {1/2}, \rlap{\qquad$\alpha\to\infty$,}\]
whereas we saw in Section~\ref{sec:example} that $s_2/s_1 \to 0$ as 
$\alpha \to\infty$ for this example.  Thus, while the results of this 
section are a marked improvement over previously existing bounds
in some highly nonnormal regimes, they cannot be the last word on the subject.
\end{remark}

\section{Conclusions}

We have illustrated a regime of stable matrices $\BA$
for which all previous bounds on the decay of singular values of
Lyapunov solutions fail to even qualitatively capture the correct
behavior.  This shortcoming is clear from specific examples;
Theorem~\ref{thm:genbnd} and Corollary~\ref{cora:genbnd} provide 
contrasting perspectives on this phenomenon.
While these results are not entirely sharp, they clearly 
illustrate that, beyond a threshold, an increased departure 
of $\BA$ from normality can lead to faster decay of the 
singular values of $\BX$.
Sharper results will require a more complete understanding 
of the role of nonnormal coefficients on Lyapunov solutions.
\new{
\section*{Acknowledgments}  
We are grateful to several referees for numerous
helpful comments on an earlier version of this
manuscript.}

\bibliographystyle{siam}
\bibliography{lyapunov}
\end{document}